\documentclass[11pt]{amsart}

\usepackage{amsmath,amssymb,amscd,verbatim} 
 \usepackage{tikz}
    \usepackage{epsfig}
    \usepackage{graphicx}
    \usepackage{color}
    \usepackage{soul}
\newlength{\guillotine}
\newtheorem{prop}[equation]{Proposition}
\newtheorem{thm}[equation]{Theorem}

\newtheorem{lemma}[equation]{Lemma}

\theoremstyle{definition}

\newtheorem{definition}[equation]{Definition}

\newtheorem{remark}[equation]{Remark}

\newtheorem{example}[equation]{Example}

\numberwithin{equation}{section}

\begin{document}
%\title{xxxx}
\date{}

\title{Zeta functions
in higher Teichm\"uller theory}

\author{Mark Pollicott and Richard Sharp}
\address{Mark Pollicott, Department of Mathematics, Zeeman Building, University of Warwick, Coventry, CV4 7AL}
\email{M.Pollicott@warwick.ac.uk}
\address{Richard Sharp, Department of Mathematics, Zeeman Building, University of Warwick, Coventry, CV4 7AL}
\email{R.J.Sharp@warwick.ac.uk}

%\affil{Warwick University}

\maketitle

\begin{abstract}
In this note we introduce 
zeta functions and $L$-functions for discrete and 
faithful representations of surface groups
in $\mathrm{PSL}(d, \mathbb R)$, for $d \geq 3$. These are natural  generalizations of  the well
known classical
 Selberg zeta function and $L$-function for Fuchsian groups, 
 corresponding to the case $d=2$.
 We  show that these complex  functions have meromorphic extensions to
 the entire complex plane  $\mathbb C$.
\end{abstract}

\section{Introduction}\label{sec:introduction}
Zeta functions and $L$-functions play a central  role in number theory, geometry and dynamical systems.  An important example is the Selberg zeta function associated to a compact hyperbolic surface.
A hyperbolic metric on a compact surface is determined by a representation of its fundamental group
into the rank-$1$ Lie group $\mathrm{PSL}(2,\mathbb R)$.
In this article, we address the natural question of how to extend  the above approach to 
zeta functions $Z(s,\rho)$ naturally  associated to appropriate representations
$\rho$ into higher rank Lie groups, 
which are the basis of the 
relatively new area of higher Teichm\"uller theory.
Here, when defining the zeta function, we replace the length of closed geodesics with the
logarithm of the spectral radius of the image of (conjugacy classes of) group elements under the representation.
A suitable  class of representations to study are those known as projective Anosov (see section 2). 

A simplified form  of one of our main results  is the following.

\begin{thm}
Let $\Gamma$ be the fundamental group
of a compact oriented surface of genus at least $2$ and let
$\rho: \Gamma \to  \mathrm{SL}(d, \mathbb R)$ ($d \geq 3$) be a projective Anosov representation. 
Then
$
Z(s, \rho)
$
converges for $Re(s)$ sufficiently large, and extends to a meromorphic function in the entire complex plane.
\end{thm}

This appears in a more precise form as  Theorem \ref{main},  which is the 
analogue of the classical result that  Selberg  zeta functions  have meromorphic extensions to the entire complex plane.

The original  Selberg zeta function may also be considered dynamically as a zeta function associated to
the geodesic flow over the surface.
 In this context, it is common to also consider the closely 
related Ruelle zeta function. Both these 
functions may be studied via a body of tools and techniques known
as Thermodynamic Formalism. In particular, this leads to results about extension beyond their half-plane of
convergence.

We now briefly outline the contents of the rest of the paper. 
In the next section, we recall the Selberg and Ruelle zeta functions and introduce analogous  
zeta functions associated to projective Anosov representations.
In section \ref{sec:surface_groups}, we discuss
surface groups and their symbolic coding. We associate to this a dynamical system, namely a subshift of finite type and explain 
how the induced action on projective space allows us 
to introduce a function defined on the shift that encodes the weightings $d_\rho(g)$.
In section \ref{extension}, we use an approach due to Ruelle \cite{ruelle},
 to obtain results on extending our zeta functions as meromorphic functions in the entire 
complex plane by writing them in terms of determinants of transfer operators. In section 
\ref{zeros}, we discuss the location of zeros and poles for the zeta functions and obtain an error 
estimate on the associated counting function. This is based on the ideas originally introduced by Dolgopyat
\cite{dolgopyat} for estimating iterates of transfer operators. A final section considers the more general case
of $L$-functions associated to unitary representations. 
Some of the technical results we require  in the paper (Lemma \ref{stuff}, Lemma \ref{key} and Lemma
\ref{mixingcondition}) originally appeared in the preprint \cite{PS}, and we include their proofs for completeness.

We would like to take this opportunity to thank the anonymous referee for their careful reading 
of the paper and many constructive comments. In particularly, we are grateful to them for 
correcting an error in the original version of section \ref{zeros}.

\section{Representations and zeta functions}\label{sec:representations}

\subsection{Selberg and Ruelle zeta functions for $\mathrm{PSL}(2, \mathbb R)$}
In 1956, Selberg introduced a zeta function associated to the fundamental group $\Gamma$ of a compact  oriented surface $V$
of genus $g \geq 2$
or, more precisely, to representations of such groups as Fuchsian groups, i.e.
discrete subgroups of $\mathrm{PSL}(2,\mathbb R)$.

The surface $V$ admits a family of hyperbolic metrics (metrics of constant curvature $-1$)
which are parametrised by the Teichm\"uller space $\mathcal T(V)$.
Let us consider  a particularly convenient viewpoint using representations.
We recall that  
$\mathcal T(V)$ can be identified with one of the two connected components of 
$$\mathrm{Hom}(\Gamma,\mathrm{PSL}(2, \mathbb R))/\mathrm{PSL}(2, \mathbb R)$$
 which consist of discrete and faithful representations
$\rho: \Gamma \to \mathrm{PSL}(2, \mathbb R)$.
(Here  $\mathrm{PSL}(2,\mathbb R)$ acts on $\mathrm{Hom}(\Gamma,\mathrm{PSL}(2, \mathbb R))$
by conjugation.)
If $\rho \in \mathcal T(V)$ then one recovers the hyperbolic metric by realising $V$
as $ \mathbb H^2/\rho(\Gamma)$, where  $\mathbb H^2$ is the Poincar\'e upper half plane
and where
$\mathrm{PSL}(2, \mathbb R)$ acts by M\"obius transformations as
orientation preserving isometries of 
$\mathbb H^2$.
(We also note that, since $\Gamma$ is torsion free, $\rho : \Gamma \to \mathrm{PSL}(2,\mathrm R)$
can be lifted to a representation $\widetilde \rho : \Gamma \to \mathrm{SL}(2,\mathbb R)$ \cite{culler}.
To simplify notation, we will use $\rho(g)$ to denote the matrix in $\mathrm{SL}(2,\mathbb R)$ 
in the hope that the context will make the usage clear.)
Another connected component $\mathcal T'(V)$ corresponds to reversing the orientation.

The Selberg zeta function  is a function of a single complex variable $s\in \mathbb C$, 
formally defined by
$$
S(s, \rho) = \prod_{n=0}^\infty \prod_{[g] \in \mathcal P} \left( 1 - e^{-(s+n) \ell_\rho(g)}\right),  \eqno(1.1)
$$
where the second product is over the set $\mathcal P$ of all primitive conjugacy classes $[g]$ of elements $g \in \Gamma \setminus \{1\}$.(The Selberg zeta function is more
usually denoted by $Z$ but we wish to reserve this for a related but slightly different function defined below.)
A conjugacy class is called primitive if it does not contain
an element of the form $g_0^n$, for $g_0 \in \Gamma$ and $n >1$.  The real number $ \ell_\rho(g)$ 
is twice the logarithm of the largest eigenvalue for $\rho(g)$ (or, equivalently, 
$\ell_\rho(g)$ is the length of the unique closed
geodesic on $V$ in the free homotopy class corresponding to $[g]$).   
The product converges to a non-zero analytic function 
for $\mathrm{Re}(s) > 1$ and Selberg showed the following  fundamental result.

\begin{thm}[Selberg]\label{selberg}
The Selberg zeta function 
$S(s,\rho)$
extends to an entire function of order $2$
 and
has a simple zero at $s=1$.
\end{thm}

An account of this theorem is contained in Hejhal's book \cite{hejhal}, where it appears as 
Theorem 4.11 and Theorem 4.25.
The result was obtained independently by Randol \cite{randol}.

One may also, following Ruelle (cf. \cite{ruelle}), consider the related zeta function
\[
R(s,\rho) = \prod_{[g] \in \mathcal P} \left(1-e^{-s\ell_\rho(g)}\right)^{-1}.
\eqno(1.2)
\] 
Since $R(s,\rho) = S(s+1,\rho)/S(s,\rho)$, one immediately obtains the following result
from Theorem \ref{selberg}.

\begin{thm}
\label{ruelle-zeta}
The Ruelle zeta function 
$R(s,\rho)$ has a meromorphic extension to $\mathbb C$ with a simple pole at $s=1$. 
Moreover, $R(s,\rho)$ can be written as a ratio of two entire functions of order $2$. 
\end{thm}

In \cite{ruelle}, Ruelle gave an alternative proof of Theorem \ref{ruelle-zeta}
which will be the basis of our approach when we  define and study analogues of these zeta functions in higher 
rank Teichm\"uller spaces.

We now want to consider a natural generalization of these definitions and results.

\subsection{Higher rank Teichm\"uller theory}
In recent years there has been considerable interest in generalising results in classical Teichm\"uller theory to
 what is now often referred to as higher Teichm\"uller theory, involving representations of surface groups
 (and more general groups) in higher rank Lie groups.    
This point of view, which had its origins in the work of Goldman \cite{Goldman88} and \cite{Hitchin87}, has  
received considerable attention, see, for example, the surveys
\cite{BCSsurvey}, \cite{kassel}
and \cite{weinhard}.  
In this note we will address the natural problem of studying analogues of the Selberg and Ruelle zeta functions in the context of higher Teichm\"uller theory.

We wish to consider  representations 
of $\Gamma$ in the higher rank group $\mathrm{PSL}(d, \mathbb R)$ (for $d \geq 3$).
The most elementary class of such representations are the so-called Fuchsian representations, obtained directly
from representations into $\mathrm{PSL}(2,\mathbb R)$.

\begin{example}[Fuchsian representations]
It is well known that there is an irreducible 
representation $\widetilde \tau: \mathrm{SL}(2, \mathbb R) \to \mathrm{SL}(d, \mathbb R)$, unique up to conjugation.
This has an explicit construction, which we briefly recall.  Let $\mathcal S_d$ denote the $d$-dimensional vector space of homogeneous polynomials in $2$ variables  of degree $\le d-1$.   We can choose a basis for $\mathcal S_d$ of the form 
$$x^{d-1}, x^{d-2}y, \cdots, x^{d-i-1} y^i, \cdots, x y^{d-2}, y^{d-1}.$$
For 
$$
g = 
\left( \begin{matrix}
a & b \\
c&d
\end{matrix}
\right) \in \mathrm{SL}(2, \mathbb R),
$$
we can define $\widetilde \tau(g) \in \mathrm{SL}(d,\mathbb R)$  by 
specifying its action on 
the basis elements as follows: 
$$
\widetilde \tau(g): x^{d-i-1} y^i \mapsto 
(ax + cy)^{d-i-1} (bx + dy)^i.
$$
(It is an easy calculation to check that $\det \widetilde \tau(g) =1$.)
Furthermore, $\widetilde \tau$ factors to an irreducible representation
$\tau: \mathrm{PSL}(2, \mathbb R) \to \mathrm{PSL}(d, \mathbb R)$.
 Representations $\rho : \Gamma \to \mathrm{PSL}(d,\mathbb R)$ of the special  form 
 $\rho = \tau \circ \rho_0$
where $\rho_0 \in \mathcal T(V)$ (or $\mathcal T'(V)$),
  are called 
\emph{Fuchsian representations}.  
\end{example}

We can now define the representations we wish to study.

\begin{definition}\label{def:hitchin_component}
The natural generalization of the 
Teichm\"uller space 
are connected components of
$$
\mathrm{Hom}(\Gamma, \mathrm{PSL}(d, \mathbb R))/\mathrm{PSL}(d, \mathbb R)
$$
which contain  Fuchsian representations.
These are  called 
\emph{Hitchin components}. When $d$ is odd there is a single Hitchin component but when $d$ is even there are two Hitchin components (corresponding to $\mathcal T(V)$ and $\mathcal T'(V)$ having distinct images).
\end{definition}

The study of representations in the Hitchin component has proved particularly fruitful.  
Hitchin studied the connected components as in 
Definition \ref{def:hitchin_component}, proving that they were each diffeomorphic to
$\mathbb R^{(2g-2)(d^2-1)}$, where $\Gamma$ is the fundamental group of a surface of genus $g \ge 2$
\cite{Hitchin87},
and further work of Goldman \cite{Goldman88} began to describe their role in parametrising geometric structures on surfaces. Labourie \cite{labour} showed that the Hitchin components consist entirely of discrete 
and faithful representations.
More recently the dynamics of the individual representations in the component have attracted attention (see, for example, the survey \cite{kassel}).

We note that any $\rho : \Gamma \to \mathrm{PSL}(d,\mathbb R)$ in a Hitchin component may
be lifted to a representation $\widetilde \rho : \Gamma \to \mathrm{SL}(d,\mathbb R)$.
This follows since a Fuchsian representation $\tau \circ \rho_0$ into $\mathrm{PSL}(d,\mathbb R)$
lifts to a representation
$\widetilde \tau \circ \widetilde \rho_0$ into $\mathrm{SL}(d,\mathbb R)$ and,
by Theorem 4.1 of 
 \cite{culler}, everything in the same connected component has the same lifting property.
 In what follows, we shall consider representations into $\mathrm{SL}(d,\mathbb R)$,
 which we
 shall denote by $\rho$ (but remembering that this is the lift of a representation into
 $\mathrm{PSL}(d,\mathbb R)$).

 \subsection{Anosov representations}
It is now convenient   to consider a more general setting.
Consider the case of discrete finitely generated torsion-free (word) hyperbolic group  $\Gamma$.  
Let $\Gamma_0$ be a finite  symmetric generating set and let $|\cdot|: \Gamma \to \mathbb Z^+$ be the word length with respect to $\Gamma_0$.

\begin{definition}\label{anosov}
A discrete and faithful representation $\rho: \Gamma \to \mathrm{SL}(d, \mathbb R)$  is called {\it k-Anosov} 
	(where  $1 \leq k \leq d/2$) if there exists $C, \mu > 0$ such that for each 
	$g \in \Gamma\setminus\{1\}$
	the eigenvalues 
	$\lambda_i(\rho(g))$, $1\leq i \leq d$, for $\rho(g)$ ordered by modulus 
(i.e., 	$|\lambda_1(\rho(g))| \geq |\lambda_2(\rho(g))| \geq \cdots \geq |\lambda_d(\rho(g))|$)
	satisfy
	$$
	\frac{\lambda_k(\rho(g))}{\lambda_{k+1}(\rho(g))} \geq C e^{|g|}.
	$$
\end{definition}

This particularly short  definition is  due to  Kassel and Potrie \cite{kp}
(see also \cite{kassel}) but corresponds to earlier definitions, beginning with 
Labourie \cite{labour} and including  Guichard and Wienhard  \cite{gw}.

A consequence of this definition is that there is a well behaved map from the Gromov boundary of $\Gamma$ to $\mathcal G_k( \mathbb R^{d})$, the Grassmannian of $k$-planes in $\mathbb R^k$.
Here, the Gromov boundary is defined as follows. 
A  sequence $(g_n)_{n=0}^\infty \in \Gamma^{\mathbb Z^+}$
 is a geodesic ray (in the Cayley graph of $\Gamma$) if 
 $|g_n|=n$ and 
 $g_{n+1}g_n^{-1} \in \Gamma_0$.
 (In particular, $g_0 = 1$.)
 Then the Gromov boundary of $\Gamma$, which we denote by
$\partial_\infty \Gamma$, 
is the space of equivalence classes of geodesic rays which remain a bounded distance apart. 
It can be equipped with a natural ``visual'' metric, where two points are close if one can choose representative rays which
 agree for a large initial segment (see Chapitre 7 of \cite{ghysdelaharpe}).
 In our case, where $\Gamma$ is the fundamental group of a compact
 surface of genus at least two, $\partial_\infty \Gamma$ is
 homeomorphic to the unit circle (by a quasisymmetric homeomorphism) \cite{ghysdelaharpe}.
 There is a natural action of $\Gamma$ on $\partial_\infty \Gamma$ and each $g \in \Gamma \setminus \{1\}$ has a unique attracting
 fixed point $g^+ \in \partial_\infty \Gamma$ and $\{g^+ \hbox{ : } g \in \Gamma \setminus \{1\}\}$
 is dense in $\partial_\infty \Gamma$.

We now have the following result.

\begin{lemma}\label{anosov}
	If a representation $\rho: \Gamma \to \mathrm{SL}(d, \mathbb R)$  is 
	$k$\emph{-Anosov} then 
	there exist H\"older continuous maps
	$$\xi^{(k)}: \partial_\infty \Gamma \to \mathcal G_k( \mathbb R^d)
	\hbox{ and }
	\theta^{(k)}: \partial_\infty \Gamma \to   \mathcal G_{d-k}( \mathbb R^{d})$$
	   such that for $\eta, \eta' \in \partial \Gamma_\infty$:
	  	\begin{enumerate}
	  		\item for $\eta \neq \eta'$  we have that 
			$\xi^{(k)}(\eta) \oplus \theta^{(k)}(\eta') = \mathbb R^d$ 
			   ({\it transversality}); 
		\item 
		$\xi^{(k)}(g\eta) = \rho(g)\theta^{(k)}(\eta)$ for all $g \in \Gamma \setminus \{1\}$ ({\it equivariance}); and 
		\item there exist $a>0$ and $c>0$ such that for every geodesic ray $(g_n)_{n=0}^\infty$ 
		corresponding  to $\eta \in \partial_\infty \Gamma$ we have that 
		$$
		\frac{\sigma_k(\rho(g_n))}{\sigma_{k+1}(\rho(g_n))} \ge ce^{an},
		$$
		where $\sigma_i$ denotes the $i$th singular value ($k$-dominated).
	\end{enumerate}
\end{lemma}

Parts (1) and (2) are standard and part (3) may be found in \cite{BPS}.

In this paper, we will 
deal with representations which are \begin{color}{blue}$1$\end{color}-Anosov. 
These representations are also called
projective Anosov. Noting that $\mathcal G_k( \mathbb R^d)$ is the projective space
$\mathbb R P^{d-1}$, in this case
we will only need the map  $\xi: \partial_\infty \Gamma \to \mathbb R P^{d-1}$ corresponding to $\xi^{(1)}$ in the lemma above. In this case, we have the following contraction property,
which, though not part of the original definitions in \cite{labour}, \cite{gw},  
comes from the equivalent  definition due to Kapovich, Leeb and Porti (\cite{klp}, Definition 6.45) and \cite{klp1}, or from part (3) of Lemma \ref{anosov}.

\begin{lemma} \label{contraction_property}
There exist $a>0$ and $c>0$ such that for every geodesic ray $(g_n)_{n=0}^\infty$ 
corresponding  to $\eta \in \partial_\infty \Gamma$ we have that 
		the associated  action $\widehat \rho(g_n): \mathbb R P^{d-1} \to  \mathbb R P^{d-1}$ satisfies 
		$\|D_{\xi(\eta)}\widehat \rho(g_n)\| \leq c e^{-n a}$, for $n \geq 0$.
\end{lemma}

Projective Anosov representations were originally introduced by Labourie \cite{labour}.
Examples include the representations in the Hitchin component introduced above and also the Benoist representations introduced in 
\cite{Benoist}.

\begin{remark}
(i) By consider the $k$th-exterior power, a \begin{color}{blue}$k$\end{color}-Anosov representation into
$\mathrm{SL}(d,\mathbb R)$ induces a projective Anosov representation into
$\mathrm{SL}\left(\left(\begin{smallmatrix} d \\k \end{smallmatrix}\right),\mathbb R
\right)$.

\noindent
(ii)
There may be additional restrictions on Anosov representations of hyperbolic groups. For example, if $\Gamma$
is a torsion-free hyperbolic group which has an Anosov representation into $\mathrm{SL}(3,\mathbb R)$
then $\Gamma$ is either a free group or a surface group (\cite{CT}, Theorem 1.1).
\end{remark}

\subsection{Zeta functions for projective Anosov representations}
We now come to the main object of study  in this note.
Let $\Gamma$ be the fundamental group of a compact oriented surface of genus $g \ge 2$ and
let
$\rho : \Gamma \to \mathrm{SL}(d,\mathbb R)$ be a projective Anosov representation.
We recall the following definition.

\begin{definition}
A matrix is called \emph{proximal}  if it has a unique eigenvalue which is real and strictly maximal in modulus.
\end{definition}

We then have the following, which is a consequence of the definition.

\begin{lemma} \label{anosov_implies_proximal}
 Let $\rho : \Gamma  \to \mathrm{SL}(d, \mathbb R)$ be 
 a projective Anosov representation. Then
   for each $g \in \Gamma \setminus \{1\}$
 the matrix $\rho(g)$ is proximal.
 \end{lemma}

We use this to associate a natural weight to each $g \in \Gamma \setminus \{1\}$.
Let   $$d_{\rho}(g) = \log |\lambda(g)|,$$ 
where $\lambda(g)$ is 
 the maximal eigenvalue of the matrix $\rho(g)$.
By Lemma \ref{anosov_implies_proximal},
since $\rho(g) \in \mathrm{SL}(d,\mathbb R)$,
$|\lambda(g)| >1$ and hence $d_\rho(g) >0$.
Note that $d_{\rho}(g)$ only depends on the conjugacy class $[g]$.
We will generalize the definition of the Selberg and Ruelle  zeta functions  to these representations.

\begin{definition}
Given a projective Anosov representation $\rho: \Gamma \to  \mathrm{SL}(d, \mathbb R)$, 
we can associate 
 a {\it generalized Ruelle zeta function} defined by
\[
\zeta(s,\rho) = \prod_{[g] \in \mathcal P} \left(1-e^{-sd_\rho(g)}\right)^{-1},
\eqno(1.2)
\]
and 
a {\it generalized Selberg zeta function} defined by 
\[
Z(s, \rho) = \prod_{n=0}^\infty \prod_{[g] \in \mathcal P} \left(
1- e^{-(s+n) d_{\rho}(g)}
\right), \eqno(1.3)
\]
for $s \in \mathbb C$ wherever the products converge.
\end{definition}

In particular, from (1.2) and (1.3) we have the identities
$$
\zeta(s,\rho) = \frac{Z(s+1,\rho)}{Z(s,\rho)}.
\eqno(1.4)
$$

\begin{remark}
Note that, in the special case $d=2$, with $\rho \in \mathcal T(V)$ we have that 
$d_\rho(g) = \ell_\rho(g)/2$. Therefore, 
$\zeta(s,\rho) = R(s/2,\rho)$.
However, the relationship between $Z(s,\rho)$ and $S(s,\rho)$
is more complicated with 
$$Z(s,\rho) = S(s/2,\rho) S(s/2 + 1/2,\rho).$$
\end{remark}

To describe the half-plane of convergence of these new zeta functions, we recall  the following definition.
See \cite{BCSsurvey}, \cite{BCLS} or \cite{PS} for a further discussion, including the existence of the limit.

\begin{definition} \label{defofentropy}
The \emph{entropy}  $h(\rho)$ of a projective Anosov representation 
$\rho : \Gamma \to  \mathrm{SL}(d, \mathbb R)$ 
is defined
to be the growth rate of the 
number of primitive conjugacy classes $[g] \in \mathcal P$
 with 
$ d_{\rho}(g)$ at most $T$ as $T \to +\infty$, i.e., 
$$
h(\rho) = 
 \lim_{T \to +\infty}  \frac{1}{T}
\log \left(
\#\{[g] \in \mathcal P \hbox{ : } d_{\rho}(g) \leq T  \}\right). \eqno(1.5)
$$
(The value of $h(\rho)$ would be the same if we used all non-trivial conjugacy classes in the definition, rather than just the primitive ones.)
\end{definition}

It is easy to  see from the definitions that $Z(s, \rho)$ and $\zeta(s,\rho)$ converge to non-zero analytic functions provided $\mathrm{Re}(s) > h({\rho})$. 
In the original setting, with $d=2$ and $\rho \in \mathcal T(V)$, the entropy is always equal to $2$ and so does
not need to be explicitly introduced. (We remark that, in this case, $h(\rho)=2$ is twice the topological entropy of the geodesic flow over $\mathbb H^2/\rho(\Gamma)$ because $d_\rho(g) = \ell_\rho(g)/2$, where 
$\ell_\rho(g)$ is the length of the corresponding closed geodesic.)

\subsection{Meromorphic extensions of the zeta functions}
Our main results on the zeta functions for projective Anosov
representations 
can be viewed as generalizations of Theorem \ref{ruelle-zeta} and Theorem \ref{selberg}, respectively.
We begin with the result for $\zeta(s, \rho)$, since, with our approach, the proof of this
 naturally comes first.

\begin{thm}\label{cor-main}
Let  $\rho: \Gamma \to  \mathrm{SL}(d, \mathbb R)$ be a projective Anosov representation. 
Then
$\zeta(s,\rho)$ extends to a meromorphic function in the entire complex plane 
and has
a simple pole at $s=h(\rho)$.
Furthermore, $\zeta(s,\rho)$ may be written as the quotient of two entire functions of order at most
$d$.
\end{thm}

 This will allow us to deduce the following result for $Z(s, \rho)$. 
 
\begin{thm} \label{main}
Let  $\rho: \Gamma \to  \mathrm{SL}(d, \mathbb R)$ be a projective Anosov representation.
Then $Z(s,\rho)$ is a meromorphic  function to the entire complex plane
and has
a simple zero  at $s=h(\rho)$.
Furthermore, $Z(s,\rho)$ may be written as the quotient of two entire functions of order at most
$d+1$.
\end{thm}

We will prove these results in section \ref{extension}.  
We will also consider the zeros of the zeta functions in section \ref{zeros}.
This leads to estimates on the error terms in the counting functions for the weights
for the subclass of hyperconvex representations.

\begin{remark}\label{rem:adjoint_rep}
An alternative way of weighting the matrices,
and hence of producing zeta functions, would be to take the difference between the logarithms of largest and smallest eigenvalues, 
i.e.
\[
 \bar {d}_\rho(g) 
=  \log\left|\lambda(g)/\lambda(g^{-1})\right|  
= \log \left(\lambda_d/\lambda_1\right).
\]
This may be reduced to the situation we consider by
composing $\rho$ with the Adjoint representation of $\mathrm{SL}(d,\mathbb R)$ and so Theorem
\ref{cor-main} and Theorem \ref{main} will still apply.
(We are grateful to the referee for this observation, which replaces our original more complicated argument.)
\end{remark}

\section{Surface groups and symbolic dynamics}\label{sec:surface_groups}
In this section we recall  the ideas that are used to build a  bridge between the geometry and the dynamics.

 \bigskip
 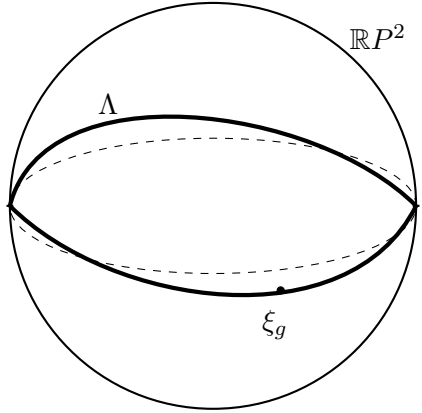
\begin{figure}[h]
   \centerline{
\begin{tikzpicture}[scale=0.45]
\filldraw (2,-2.5) circle (3pt);
\draw[thick] (0,0) circle (6cm);
\node [left] at (6,5) {$\mathbb R P^2$};
\node [left] at (2.5,-3.5) {$\xi_g$};
\node [left] at (-2.5,3) {$\Lambda$};
   \draw[thick] (5.9,0) --  (6.1,0);
      \draw[thick] (-5.9,0) --  (-6.1,0);
      \draw[dashed] (0,0) ellipse (6cm and 2cm);
      \draw[ultra thick] (-6,0) .. controls (-3, -3) and (4, -4) .. (6, 0);
            \draw[ultra thick] (-6,0) .. controls (-5, 4) and (3, 3) .. (6, 0);
\end{tikzpicture}
}
\caption{
Here $\rho$ is a representation into
$\mathrm{SL}(3,\mathbb R)$ and $\Lambda  \subset \mathbb R P^2$ is the
associated limit set.  For $g \in \Gamma \setminus \{1\}$, the point 
$\xi_g$ is an attracting fixed point for $\widehat \rho(g): \mathbb R P^2 \to : \mathbb R P^2$.}
  \end{figure}

 \subsection{Action on projective space}
 As discussed above, the representation $\rho : \Gamma \to \mathrm{SL}(d,\mathbb R)$ induces an action of $\Gamma$
 on the projective space $\mathbb R P^{d-1}$. More explicitly, for $g \in \Gamma$, we define
 $\widehat \rho(g): \mathbb R P^{d-1}  \to \mathbb R P^{d-1}$ 
by 
 \[
 \widehat \rho(g)[v] = \frac{\rho(g)v}{\|\rho(g)v\|_2},
 \]
 where $v \in \mathbb R^d \setminus \{0\}$ is a representative element.
 
 Suppose $g \in \Gamma \setminus \{1\}$.
 Since, by Lemma \ref{anosov_implies_proximal}, $\rho(g)$ has a real eigenvalue which is strictly maximal in modulus,
 the map $\widehat \rho(g)$ has a unique attracting fixed point; we will denote this by
$\xi_g \in \mathbb R P^{d-1}$.
Let $\xi : \partial_\infty \Gamma \to \mathbb R P^{d-1}$ be the map introduced in the previous section.
 As an immediate consequence of Lemma \ref{contraction_property}, for each 
 $g \in \Gamma \setminus \{1\}$,
we have
 $\xi(g^+) = \xi_g$, the fixed point for $\widehat \rho(g)$.

We now define a compact subset of projective space which will be useful in the sequel.

\begin{definition}
We define the {\it limit set} 
$$\Lambda:= \overline {\{\xi_g \hbox{ : } g \in \Gamma \setminus \{1\}\}} \subset \mathbb R P^{d-1}$$
to be the closure of the fixed points.
\end{definition}

Since $\{g^+ \hbox{ : } g \in \Gamma \setminus \{1\}\}$ is dense in $\partial_\infty \Gamma$,
immediately follows that
$\Lambda = \xi(\partial_\infty \Gamma)$.

\begin{example}[Case $d=3$]
In this case the limit set $\Lambda \subset \mathbb R P^2$ is diffeomorphic to a  circle bounding a strictly convex region
$\Omega$
 when $\rho$ a representation in the Hitchin component. 
   In particular, in Lemma \ref{anosov} we have that 
    for $\eta \in \partial_\infty \Gamma$ the point $\xi^{(1)}(\eta) \in \Lambda$ and 
$\theta^{(1)}(\eta)$ is tangent to the curve at that point.
 This is described in the earlier  work of Choi and Goldman  \cite{CG} (see also \cite{labour}) and applies to projective Anosov representations via the equivalence of the  definitions shown in  \cite{kp}.
\end{example}

We can use the following simple lemma 
to  relate the weight $d_\rho(g)$ to the action of 
$\widehat \rho(g)$ on 
$\mathbb R P^{d-1}$.
 
 \begin{lemma}  \label{stuff}
If $g \in \Gamma \setminus \{1\}$   
then 
\[
d_\rho(g) =  -\frac{1}{d} \log \det(D_{\xi_g} \widehat \rho(g)).
\]
\end{lemma}

\begin{proof}
We can consider the linear action of $\rho(g)$  on $\mathbb R^d$, then  the fixed point  
corresponds to an eigenvector $v$ and the result follows from a simple  calculation using that
the linear action of  $\rho(g) \in \mathrm{SL}(d, \mathbb R)$ preserves area in $\mathbb R^d$.
More precisely, $\xi_g$ corresponds to an eigenvector $v$ for the maximal eigenvalue $\lambda(g)$,
with $|\lambda(g)| > 1$, for the matrix 
$\rho(g)$.  We can assume without loss of generality  that $\|v\|=1$ and  then for arbitrarily small $\delta>0$ we can consider a $\delta$-neighbourhood of $v$ which is the product of a
$(d-1)$-dimensional neighbourhood in $\mathbb R P^{d-1}$ and a $\delta$-neighbourhood
in the radial direction. 
The effect of the linear action of $\rho(g)$ is to replace $v$ by $\lambda(g) v$,  and thus 
stretch the neighbourhood in the radial direction by a factor of $|\lambda(g)|$.  
Since $\rho(g)$ has determinant one, the volume of the $(d-1)$-dimensional neighbourhood
contracts by $|\lambda(g)|^{-1}$. 
To calculate the effect of the projective action $\widehat \rho(g)$, we need to rescale  
$\lambda(g) v$ to have norm one, which corresponds to multiplication by the diagonal
matrix $\mathrm{diag}(|\lambda(g)|^{-1},\ldots,|\lambda(g)|^{-1})$.
 In particular, the $(d-1)$-dimensional neighbourhood in $\mathbb R P^{d-1}$ shrinks by a factor of approximately 
 $|\lambda(g)|^{-d}$, giving the result.
\end{proof}

\subsection{Surface groups and coding}
We now want to introduce a natural coding for the group $\Gamma$ which will allow us  to 
analyse the projective action of $\Gamma$ on $\Lambda$.
This is motivated by the more familiar boundary coding associated to Fuchsian groups due to Bowen and Series \cite{series}.
We  take the same basic approach 
to introducing a dynamical perspective to the study of projective Anosov representations 
 based
on subshifts of finite type that was used in \cite{PS} (for representations in a Hitchin component).

Since $\Gamma$ is the fundamental group of a compact surface $V$ with genus $g \geq 2$, it has the presentation 
$$\Gamma 
= \left\langle a_1,\ldots,a_g, b_1,\ldots,b_g \hbox{ : }  \prod_{i=1}^g[a_i, b_i] = 1\right\rangle.$$
We write
 $\Gamma_0 =\{a_1^{\pm 1}, \cdots, a_g^{\pm 1}, b_1^{\pm 1}, \cdots, b_g^{\pm 1} \}$
 for the symmetrized generating set and $|\cdot| : \Gamma \to \mathbb Z^+$ for the word 
 length with respect to $\Gamma_0$.

Such surface groups are particular examples of 
Gromov hyperbolic groups and as such they are strongly Markov groups in the sense of Ghys and de la Harpe \cite{ghysdelaharpe}, i.e. they can be encoded using a directed graph and an edge labelling by elements in $\Gamma_0$. This approach allows allows us to relate the group to a dynamical system, namely a subshift of finite 
type, and thus to use the machinery of thermodynamic formalism to define a so called pressure 
function and, later, transfer operators, which may be used to analyse our zeta functions.
In the particular case of surface groups, the coding follows directly from the 
work of Adler and Flatto \cite{AF} and Series \cite{series} on coding the action on the boundary
and the associated subshift of finite type is mixing.

We have the following result.

\begin{lemma}
\label{hyp-sm} We can associate to $(\Gamma,\Gamma_0)$
\begin{enumerate}
\item[(i)] a directed graph $\mathcal G = (\mathcal V,\mathcal E)$ with a distinguished vertex $*$; and
\item[(ii)] an edge labelling $\omega : \mathcal E \to \Gamma_0$,
\end{enumerate}
such that 
\begin{enumerate}
\item no edge terminates at $*$;
\item there is at most one directed edge joining each ordered pair of vertices;
\item the map from the set finite paths in the graph starting at $*$ to $\Gamma \setminus \{1\}$ defined by
$$(e_1, \ldots , e_n) \mapsto \omega(e_1)\cdots \omega(e_n)$$
is a bijection and $|\omega(e_1)\cdots \omega(e_n)|=n$;
\item
the map from closed paths in $\mathcal G$ to conjugacy classes to $\Gamma$ induced by $\omega$ 
is a bijection and for such a closed path $(e_1, \ldots , e_n,e_1)$, $n$ is the  minimum word 
length in the conjugacy class of $\omega(e_1)\cdots \omega(e_n)$.
\end{enumerate}
Furthermore, the subgraph obtained be deleting the vertex $*$ has the aperiodicity property that
there exists $N \geq 1$ such that, given any two $v,v' \in \mathcal V \setminus \{*\}$, there is a directed path of length $N$ from $v$ to $v'$.
\end{lemma}

\begin{remark}
Part (4) of Lemma \ref{hyp-sm} plays a subtle, but important, role in allowing us to relate the dynamical theory of zeta functions to conjugacy classes in $\Gamma$.  In order to study other groups $\Gamma$ than surface groups one would need to find a replacement for this property.
\end{remark}

\subsection{Symbolic dynamics}
We now introduce the symbolic coding for the limit set $\Lambda$. 
 This is particularly useful in introducing a partition of $\Lambda$ which gives an effective way to reduce the group 
 action to the orbits of a single transformation.

We can associate to the directed graph $\mathcal G = (\mathcal V,\mathcal E)$ a subshift of finite type $\Sigma$, 
as follows. Let $\mathcal E^* \subset \mathcal E$ be the set of directed edges that start at $*$.
The set of states in $\Sigma$ is $\mathcal E \setminus \mathcal E^*$ and a transition
$e \to e'$ is allowed if and only if $e'$ follows $e$ in $\mathcal G$. 
More precisely, if
there are $k$ edges in $\mathcal E \setminus \mathcal E^*$ then we define a $k \times k$ matrix $A$ by $A(e,e') =1$ if 
$e'$ follows $e$ in the directed graph and define 
$$
\Sigma = \{x = (x_n) \in \{1, \cdots, k\}^{\mathbb Z^+} \hbox{ : } A(x_n, x_{n+1})=1, n \geq 0\}.
$$
(We note that this definition does not involve the edge labels on $\mathcal E$.)
The set $\Sigma$
is a compact space with respect to the metric
$$
d(x,y) = \sum_{n=0}^\infty \frac{1- \delta(x_n, y_n)}{2^n}.
$$
The shift map is the local homeomorphism $\sigma: \Sigma \to \Sigma$ defined by 
$(\sigma x)_n = x_{n+1}$. By the final statement in Lemma \ref{hyp-sm}, $A$ is aperiodic (i.e. there exists 
$N \geq 1$ such that $A^N$ has all
entries positive) and, equivalently, the shift $\sigma: \Sigma \to \Sigma$ is mixing (i.e. for all open
non-empty $U,V \subset \Sigma$, there exists $N \geq1$ such that 
$\sigma^{-n}(U) \cap V \neq \varnothing$ for all $n \geq N$).

There is a natural surjective H\"older continuous map
$\pi : \Sigma \to \partial_\infty \Gamma$ defined by letting $\pi((x_n)_{n=0}^\infty)$ be the equivalence classes of 
the geodesic ray $(g_n)_{n=0}^\infty$ in $\partial_\infty \Gamma$, 
where $g_n = \omega(x_0) \omega(x_1) \cdots \omega(x_n)$.

\begin{definition}
We define a map $\kappa : \Sigma \to \Lambda \subset \mathbb R P^{d-1}$
by $\kappa = \xi \circ \pi$.
\end{definition}

However,  the shift $\sigma: \Sigma \to \Sigma$ 
only encodes information about $\Gamma$ as an abstract group.  In order to keep track of  the  additional  information given by the representation of $\Gamma$  in $\mathrm{SL}(d, \mathbb R)$ 
and its action on $\mathbb RP^{d-1}$ we need to introduce a H\"older continuous function 
$r: \Sigma \to \mathbb R$.

 \begin{definition}
We can associate a  map  $r: \Sigma \to \mathbb R$ defined by 
\[
r(x) = -\frac{1}{d} \log \det(D_{\kappa(x)}\widehat \rho(g_{x_0}))
\]
 (i.e. the Jacobian of the derivative of the projective action), where $g_{x_0} = \omega(x_0)$ is the 
 generator corresponding to the first term in $x = (x_n)_{n=0}^\infty \in \Sigma$.
\end{definition}

There is a natural one-to-one correspondence between non-trivial conjugacy classes 
in $\Gamma$ and periodic orbits for the shift $\sigma$, by part (4) of Lemma \ref{hyp-sm}.  

We now have the following simple but key result.

\begin{lemma}  \label{key}  The  function  $r : \Sigma  \to \mathbb R$ is H\"older continuous, and if
 $\sigma^nx =x$ is a periodic point corresponding to the conjugacy class
 $[g]$ of an element $g \in \Gamma$  then $r^n(x)  = d_\rho(g)$.
\end{lemma} 

\begin{proof}
The H\"older continuity  of $r$ follows immediately from the H\"older continuity of $\xi$, 
which is given by Lemma \ref{contraction_property}.
   The second statement  follows from the equivariance and the observation 
   $\kappa(\sigma x)= \rho(g_{x_0}) \kappa(x)$.  Moreover,    
   that the periodic point $x$ has an image $\kappa(x) (= \xi_g)$ which is fixed by 
   $\widehat \rho (g)$ and the result follows from Lemma \ref{stuff}.
 \end{proof}

\subsection{Semi-flows and entropy}
In this subsection we give an alternative definition of the entropy 
$h(\rho)$ of  a representation, as defined in the introduction,  which is useful when taking  the symbolic viewpoint.

We can  associate to $\sigma: \Sigma \to \Sigma$ and the positive function $r: \Sigma \to \mathbb R^+$ the space
$$
\Sigma^r = \{(x, u) \in \Sigma \times \mathbb R \hbox{ : } 0 \leq u \leq r(x)\}/\sim,
$$
where $(x,r(x)) \sim (\sigma x, 0)$, and the semi-flow $\sigma^r_t: \Sigma^r  \to \Sigma^r$ defined by 
$\sigma^r_t(x,u) = (x, u+t)$, for $t \geq 0$ (respecting the identifications).

The topological entropy $h(\sigma^r)$ for the semi-flow is defined to be the entropy of the time one map.
It is a standard result that an alternative way of  defining $h(\sigma^r)$ is as the growth rate
\[
h(\sigma^r) = \lim_{T \to +\infty} \frac{1}{T} \log \#\{\tau \in \mathcal P(\sigma^r) \hbox{ : } l(\tau) \le T\},
\]
where $ \mathcal P(\sigma^r)$ denotes the set of prime periodic orbits of $\sigma^r$ and $l(\tau)$ denotes the least period of the orbit $\tau$ \cite{PP}.

\begin{lemma}\label{zero}
If $\rho$ is a projective Anosov representation then 
the entropy of the representation is equal to the topological entropy of the semi-flow $\sigma^f$.
\end{lemma}

\begin{proof}
We observe that there is a one-one correspondence between the periodic orbits for the semi-flow $\sigma^r$ and for the shift $\sigma$.  
By Lemma \ref{key} 
and Definition \ref{defofentropy}, we see that $h(\rho) =h(\sigma^r)$.
\end{proof}

\section{Extensions of Zeta Functions} \label{extension}
We now turn to the results on the meromorphic extensions of zeta functions.
The proof of Theorem  \ref{cor-main} 
 is based on the approach of Ruelle in 
\cite{ruelle}.
 We then use 
 Theorem  \ref{cor-main} to prove
 the corresponding result for $Z(s, \rho)$ in Theorem \ref{main}.

\subsection{Partitions}
By the construction in the previous section
we have associated to the representation $\rho$
a H\"older continuous map $\kappa: \Sigma \to \mathbb R P^{d-1}$ which intertwines the shift map 
on $\Sigma$ with the projective action, 
i.e. $\kappa(\sigma x) = \widehat \rho(\omega(x_0)) \kappa(x)$.
  In particular, we use $\kappa$ to define a particular partition $P_i = \kappa([i])$ of the limit set 
$\Lambda \subset \mathbb RP^{d-1}$, 
where 
$$[i] = \{x \in \Sigma \hbox{ : } x_0=i\}, \quad
1 \leq i \leq k.$$ 
We then have that
\begin{enumerate}
\item $\Lambda = \bigcup_{i=1}^k P_i$, and
\item $P_i \cap P_j$ consists of at most one point,
for distinct $i \neq j$.
\end{enumerate}

The analytic extension of the zeta function $\zeta(s,\rho)$ is based on the study of a transfer operator defined on analytic functions.  In particular, we can realise the symbolic maps 
$\Sigma \supset [j] \ni x \mapsto ix \in [i] \subset \Sigma$ 
(where $A(i,j)=1$)
as local real analytic contractions 
\[
\psi_i :
\coprod_{A(i,j)=1} P_j
\to P_i
\]
on the disjoint union of the sets.

 \bigskip
 \begin{figure}[h]
   \centerline{
\begin{tikzpicture}[scale=0.45]
\draw[thick] (0,0) circle (6cm);
\node [left] at (6,5) {$\mathbb R P^2$};
\node [left] at (-3.5,1.5) {$\Lambda$};
   \draw[thick] (5.9,0) --  (6.1,0);
      \draw[thick] (-5.9,0) --  (-6.1,0);
      \draw[semithick] (-6,0) .. controls (-3, -3) and (4, -4) .. (6, 0);
            \draw[semithick] (-6,0) .. controls (-5, 4) and (3, 3) .. (6, 0);
                  \draw[semithick, rotate=5] (1.5,-2.5) ellipse (3cm and 1.2cm);
                      \draw[ultra thick] (0.8,-2.6) .. controls (1.2, -2.7) and (2.7, -2.6) .. (3,-2.3);
                                          \node[left] at (-0.3,-1.2) {$U_j$};
         \draw[semithick, rotate=-20] (1.5,2.2) ellipse (3.8cm and 1.5cm);
                      \draw[ultra thick] (0.5,2.5) .. controls (1.0, 2.3) and (2.7, 2.0) .. (3.5,1.6);
                    \node[left] at (0.9,-2) {$P_j$}; 
                                            \node[left] at (0.3,3.9) {$U_i$};      
                                                     \draw[semithick, rotate=-20, fill=yellow, opacity=0.5] (2,2.6) ellipse (1.2cm and 0.5cm);           
                                                        \node [left] at (1.5,2.9) {$P_i$};    
                                                         \node [left] at (4.5,0.6) {$\psi_i(U_j)$};               
\end{tikzpicture}
}
\caption{The contractions $\psi_i : P_j \to P_i$ satisfy $\psi_i(U_j) \subset U_i$ for the neighbourhoods of the complexifications (suggested in the figure).}
  \end{figure}
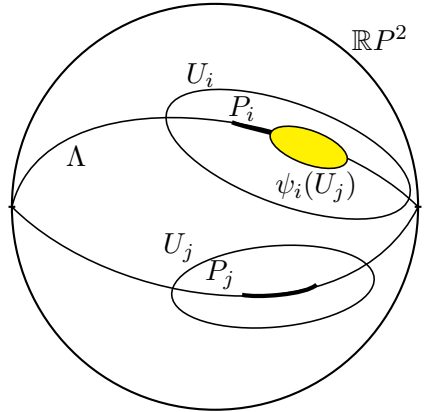

\subsection{Transfer operators}
We can consider a  Banach space of analytic functions on a neighbourhood of the  
disjoint union of the 
partition $\coprod_{i=1}^k P_i$.   In particular, 
let us fix neighbourhoods
 $U_i \supset P_i$  in the complexification of 
$\mathbb R P^{d-1}$, chosen suitably small as will be indicated below. 
 
 This allows us to define a convenient   Banach space of analytic functions.
 
\begin{definition}
We define $B_0$ to be the space of bounded analytic functions 
$f: \coprod_{i=1}^k U_i \to \mathbb C$ with the supremum norm
$$
\|f \| = \max_i \sup_{z\in U_i} |f(z)|, 
$$
i.e., $B_0 = \left\{f: \coprod_{i=1}^k U_i \to \mathbb C \ \mathrm{ analytic}\hbox{ : } \|f\| < \infty\right\}$.
\end{definition}

We can also define a natural family of bounded linear operators on this space.

\begin{definition}
For $s \in \mathbb C$, the {\it transfer operator}  
$\mathcal L_s : B_0 \to B_0$ is defined by 
$$
\mathcal L_s f(z) =  
\sum_{A(i,j)=1} \det(D\psi_i(z))^s f(\psi_iz), 
\hbox{ for } z \in U_j,
$$
where $\det (D\psi_i (\cdot ))^s$ is understood as the complexification of the real analytic function.
\end{definition}

In order for the operator to be well defined we need to show that 
we can choose the neighbourhoods $U_i$ so that
whenever $A(i,j)=1$
($1 \leq i,j \leq k$) then
 $\overline {\psi_i(U_j)} \subset U_i $.  Since, by construction, 
 $\psi(P_j) \subset P_i$, it suffices to show that for any $\eta \in P_i$ we have  $\|D_\eta \psi\| < 1$. However, associating to $\eta$ 
 those geodesics $(\gamma)$  arising from the coding (in Lemma \ref{hyp-sm}) we can apply 
 Lemma \ref{contraction_property}.  If $ae^{-c}<1$ then we are done.  Otherwise we can choose $n\geq 2$ sufficiently large so that $ae^{-cn}$ and use the standard recoding of  $\Sigma$ by replacing symbols by words of length $n$.

The advantage of being able to work in the Banach space of bounded analytic functions is that the operators have far better spectral properties than, say, Banach spaces of $C^1$ functions or H\"older continuous functions.  
More precisely, the transfer operators are compact operators and, furthermore, satisfy the stronger property of nuclearity.
 This follows directly from the work of Ruelle, based on the ideas of Grothendieck.
 
\begin{lemma}[Nuclearity: after Grothendieck and Ruelle] \label{nuclear}
The transfer operators $\mathcal L_s : B_0 \to B_0$ are nuclear.  More precisely,  we can find 
\begin{enumerate}
 \item
 $v_m \in B_0$ with $\|v_m\|_{B_0}=1$,
  \item
 $l_m \in B_0^*$ with $\|l_m\|_{B_0^*}=1$, and
 \item
 $\lambda_m \in \mathbb C$ with $|\lambda_m| = O(e^{-cm})$, for some $c>0$,
 \end{enumerate}
 such that
 \[
 \mathcal L_s f= \sum_{m=1}^\infty \lambda_m l_m(f) v_m.
 \]
 \end{lemma}

 \begin{proof}
 This can be deduced from \cite{ruelle}.
 \end{proof}
  
We can also consider the natural generalisation of these operators to the space of analytic
$j$-forms on $\coprod_{i=1}^k U_i$ (for $j=1,\ldots,  d-1$). More precisely, let 
$B_j$ denote the Banach space of analytic $j$-forms on $\coprod_{i=1}^k U_i$ whose coefficients
are bounded.
Then the definition of $\mathcal L_s$ naturally extends to define transfer operators
$\mathcal L_{s,j} : B_j \to B_j$. The analogue of Lemma \ref{nuclear} holds for these operators.

 \subsection{Traces and determinants}
A particular
consequence of Lemma \ref{nuclear} is that the operator $\mathcal L_s$ (and its powers  $\mathcal L_s^n$) are trace class, i.e. the eigenvalues are summable:
$\sum_{m=1}^\infty |\lambda_m| < \infty$.

To calculate the trace of $\mathcal L_s^n$, we need to consider compositions of the contractions $\psi_i$ corresponding to periodic points in $\Sigma$. 
Let $\underline i = (i_1,\ldots,i_n)$ be a string of length $|\underline i|=n$ an allowed string such 
that
$A(i_1, i_2) = \cdots = A(i_{n-1}, i_n)=1$. If, in addition, $A(i_n,i_1) = 1$
then we call $\underline i$ a periodic allowed string. We let 
\[
\psi_{\underline i} : \coprod_{A(i_n,j)=1} P_j \to P_{i_1}
\]
denote the composition $\psi_{\underline i} = \psi_{i_1} \circ \cdots \circ \psi_{i_n}$.
Since $A(i_n,i_1)=1$, this restricts to a contraction 
$\psi_{\underline i} : P_{i_1} \to P_{i_1}$, which has a unique fixed point denoted $z_{\underline i}$.
We now have the following result.

   \begin{lemma}\label{traces}
 The trace of $\mathcal L_s^n: B_0 \to B_0$ is given by
 $$
 \hbox{\rm tr}(\mathcal L_s^n)
 = 
 \sum_{|\underline i|=n} 
 \frac{\det (D\psi_{\underline i}(z_{\underline i}))^s}{ \det (I - D\psi_{\underline i}(z_{\underline i}))},
 $$
 where the summation is taken over all periodic allowed strings of length $n$.
 Furthermore, 
$$
\hbox{\rm tr}(\mathcal L_{s}^n) + 
\sum_{j=1}^d
(-1)^j \hbox{\rm tr}(\mathcal L_{s,j}^n)
 = 
 \sum_{|\underline i|=n} 
 \det (D\psi_{\underline i}(z_{\underline i}))^s.   \eqno(4.1)
$$
  \end{lemma}
  
\begin{proof}
 This can be deduced from the general setting in \cite{ruelle},
 pages 232--234.
 \end{proof}
 
 We want to associate to the transfer operators a function of a single complex variable $s$.
  We can write
  $$
  \det(I -  \mathcal L_s) = \exp \left( - \sum_{n=1}^\infty  \frac{1}{n}  
   \hbox{\rm tr}(\mathcal L_s^n)
   \right).
  $$ 
  Lemma \ref{traces} explains the connection  with the fixed points of the contractions.
  Similarly,  we can associate to the transfer operators $\mathcal L_{s,j}$ the functions
    $$
  \det(I -  \mathcal L_{s,j}) = \exp \left( - \sum_{n=1}^\infty  \frac{1}{n}  
   \hbox{\rm tr}(\mathcal L_{s,j}^n)
   \right)
  $$ 
 for $1 \leq j \leq d-1$.        
  
The  next lemma shows that these  functions  extend analytically to the entire complex plane.
To simplify notation, we shall write $\mathcal L_{s,0}=\mathcal L_{s}$.

\begin{lemma} \label{det}
The functions 
$s  \mapsto \det(I - \mathcal L_{s,j})$ ( $0 \leq j \leq d-1$) 
extend  to entire analytic functions of $s\in \mathbb C$ of order at most
$d$.
\end{lemma}
\begin{proof}
 This can be (essentially)  deduced from 
 the proof of Theorem 1(b) in \cite{ruelle},
 although there is a minor correction required, as described in 
 \cite{fried}.
 In particular, one can expand 
 $$
  \det(I -  \mathcal L_{s,j}) =
  1 + 
  \sum_{m=1}^\infty c_m(s) 
 $$
 where 
 \[
 c_m(s) = 
 O \left(  
 e^{|s| Dm} 
 e^{
 -Em^{\left(1+ \frac{1}{d-1} \right)}
 }
 \right),
 \] 
 for some $E>0$
 and $D = \max_{j}\|D\phi_j\|$.  
 In particular, by considering critical points of 
the real function 
$t\mapsto |s| D t
 -E t^{\left(1+ \frac{1}{d-1} \right)}$ 
 we see that the upper bound on the terms occurs at 
$$m 
 = \left[ 
 \left(
 \frac{|s| D}{E(1 + \frac{1}{d-1})}
 \right)^{d-1}
 \right]
 \hbox{ or }
 \left[ 
 \left(
 \frac{|s| D}{E(1 + \frac{1}{d-1})}
 \right)^{d-1}
 \right] + 1.
 $$
 giving an upper bound on the terms of 
 $O\left(
  \exp 
  \left(
  |s|^d D^d/(E (1 + \frac{1}{d-1}))^{d-1}
  \right)
  \right)$.  
 Moreover, it is easy to see that this gives an upper bound 
 \[
 |\det(I - \mathcal L_{s,j})| = O\left(e^{c|s|^{d}}\right),
 \] 
 for some $c>0$.
In particular, we see that the order  is given by
 $$
 \limsup_{R \to +\infty}  \frac{\log \left(\max_{|s| \leq R} \log |\det(I-\mathcal L_{s,j})| \right)}{\log R}
 \leq  d.
 $$
 \end{proof}

\noindent
{\it Proof of Theorem \ref{cor-main}.}
To complete the proof of Theorem \ref{cor-main},  we can then use the identity (4.1) to write
   $$
  \zeta(s, \rho) =  
  \frac{\prod_{j \enskip \mathrm{odd} }\det(I- \mathcal L_{s,j})}{
  \prod_{j \enskip \mathrm{even}}
\det(I- \mathcal L_{s,j}) 
  },
  $$
  from which we deduce that $\zeta(s, \rho)$ is a meromophic function which, by Lemma \ref{det}, may be written as the quotient of 
  two entire functions of order $\le d$.

\medskip
\noindent
{\it Proof of Theorem \ref{main}.}
 Given $s_0 \in \mathbb C$, we can choose 
 $m \in \mathbb Z$ with $\mathrm{Re}(s_0) \geq m$.  
Setting $n = |m| + [h(\rho)] + 1$ 
we can use (1.4) to write 
$$
Z(s, \rho) = \frac{Z(s+n, \rho)}{\zeta(s, \rho)
\zeta(s+1, \rho) \cdots \zeta(s+n-1, \rho)}.
\eqno(4.2)
$$
For $s$ in a sufficiently small
 neighbourhood of $s_0$ we have that $\mathrm{Re}(s+n) > h(\rho)$ and therefore  $Z(s+n, \rho)$ converges to an analytic function.  Moreover, the denominator is meromorphic by Theorem \ref{cor-main} and thus using (4.2) we deduce that $Z(s, \rho)$ is meromorphic.

To obtain the order  bound in Theorem \ref{main}
consider 
$|s| < m$, say, and set $n= m + [h(\rho)]+2$. 
By the proofs of Theorem \ref{cor-main}
and Lemma \ref{det}
we can write
$\zeta(s, \rho)= \frac{f(s)}{g(s)}$,  where 
for any $\epsilon>0$ we can choose $C > 0$ such that
$|f(s)|, |g(s)| \leq \exp(C |s|^{d+\epsilon})$ and thus using (4.2) we can write
$$ Z(s, \rho)= Z(s+n, \rho)
\frac{G(s)}{F(s)}\eqno(4.3)
$$
where $F(s) = \prod_{k=0}^{n-1}f(s+k)$ and $G(s) = \prod_{k=0}^{n-1}g(s+k)$ are entire functions.   In particular,
$Z(s+n, \rho)$ is uniformly bounded and since 
$$
\log |F(s)| \leq C \sum_{k=0}^{m+1} |s+k|^{d+\epsilon}
\leq C (m+2) (2m)^{d+\epsilon}
$$
the order of $F(s)$ is bounded by 
$$\lim_{m\to +\infty}\frac{\log(C (m+2) (2m)^{d+\epsilon})}{\log m} = d+1+ \epsilon,$$  and similarly the order for $G(s)$ is bounded by $d+1 +\epsilon$.  Finally, since $\epsilon > 0$ can be arbitrarily 
small, we can   deduce from (4.3) that 
$Z(s, \rho)$ is a ratio of two entire functions
of order $d+1$.

 \begin{remark}[Equidistribution]
 For each $t > \delta$ can define a probability measure on $\Lambda$ by 
 $$
\mu_t = 
\frac{
\sum_{g \in \Gamma \setminus \{1\}}  |\hbox{\rm Jac} D_{\xi_g} \widehat \rho (g)|^t \delta_{\xi_g}
}{
\sum_{g \in \Gamma \setminus \{1\}}  |\hbox{\rm Jac} D_{\xi_g} \widehat \rho (g)|^t 
}
$$
and then 
as $t$ tends to $\delta$ this converges to a probability measure $\mu$ supported on $\Lambda$.  
This follows easily from the properties of the transfer operator $\mathcal L_s$, where 
$\mu$ is related to the maximal eigenprojection of the operator $\mathcal L_{h(\rho)}$.
 \end{remark}
 
\section{Zeros and poles of the zeta functions} \label{zeros}

Having established that 
the zeta functions 
$\zeta(s, \rho)$ and 
 $Z(s, \rho)$ have  meromorphic  extensions to $\mathbb C$ in  Theorems  \ref{cor-main}
 and \ref{main}, it is very natural  to ask about the location of the zeros 
 and poles of zeta function.

 In the case of the classical Selberg zeta function
for hyperbolic  surfaces the following is well known. 
Recall from Theorem \ref{selberg} that the Selberg zeta function for a hyperbolic surface is entire with a simple zero at $s=0$.  Furthermore, the following result is well known.
  
  \begin{prop}[\cite{hejhal}]\label{riemannhypothesisselberg} 
There exists 
$\epsilon > 0$ such that $S(s,\rho)$ has no zeros in 
$\mathrm{Re}(s) >  1 - \epsilon$  other than $s=1$. 
\end{prop}

In this classical case it is even known that the zeros in the critical strip $0 \leq \mathrm{Re}(s) \leq 1$ lie on $[0,1] \cup \left(\frac{1}{2}+i \mathbb R\right)$.  This is a consequence of the way the zeta function can be extended using the trace formula and the interpretation of the zeros
$s_n$ in terms of the 
eigenvalues $\lambda_n = s_n(1-s_n)$ 
 of the Laplace-Beltrami operator.

Using the relation $R(s, \rho) = S(s+1, \rho)/S(s,\rho)$ we can deduce the corresponding 
result for $R(s,\rho)$.

  \begin{prop}\label{riemannhypothesisruelle} 
There exists 
$\epsilon > 0$ such that $R(s,\rho)$ has no 
zeros or poles in 
$\mathrm{Re}(s) >  1 - \epsilon$  other than 
a simple pole $s=1$. 
\end{prop}

We can  generalize both Proposition 
\ref{riemannhypothesisselberg} and 
Proposition \ref{riemannhypothesisruelle} to the zeta functions for 
a special class
of projective Anosov representations.
We say that $\rho : \Gamma \to \mathrm{SL}(d,\mathbb R)$ is $(1,2)$-Anosov if it is $1$-Anosov
and $2$-Anosov. A $(1,2)$-Anosov representation is $(1,1,2)$-\emph{hyperconvex} if
for every triple of pairwise distinct points $x,y,z \in \partial_\infty\Gamma$ we have
\[
\xi^1(x) \oplus \xi^1(y) \oplus \xi^{d-2}(z) = \mathbb R^d.
\]
This is an open condition and examples include representations in Hitchin components \cite{labour}.
Importantly, for our analysis,
recent work of 
Pozzetti, Sambarino and Weinhard
(Proposition 7.4 and Corollary 7.6
of \cite{PSW})
and Zhang and Zimmer (Theorem 1.1 of \cite{ZZ}) 
showed that if $\rho$ is a $(1,1,2)$-hyperconvex representation then
$\xi^1(\partial_\infty \Gamma)$ is a $C^1$ curve. Indeed, Zhang and Zimmer prove stronger 
regularity, which we state this as a lemma.

\begin{lemma}[Zhang and Zimmer, \cite{ZZ}]\label{lem:limit_set_is_C1}
If $\rho : \Gamma \to \mathrm{SL}(d,\mathbb R)$ is
$(1,1,2)$-\emph{hyperconvex} then the limit set
$\Lambda = \xi^1(\partial_\infty \Gamma)$ is a $C^{1+\alpha}$ curve, for some $\alpha>0$.
\end{lemma}

Another ingredient is the following ``mixing condition'' on the weights $d_\rho(g)$.

\begin{lemma}[Carvajales, \cite{Carvajales}, Lemma A.2]\label{mixingcondition}
Let $\rho$ be a projective Anosov representation.
Then there does not exist $a > 0$ such that $\{d_\rho(g) \hbox{ : } g \in \Gamma \setminus \{1\}\} \subset a \mathbb Z$.
\end{lemma}

We will now obtain our results on the zeta functions.
First we will consider the Ruelle-type zeta function
$\zeta(s, \rho)$.

\begin{thm}\label{riemannhypothesis} 
Let $\rho : \Gamma \to \mathrm{SL}(d,\mathbb R)$ be a $(1,2)$-Anosov representation
which is $(1,1,2)$-hyperconvex.
Then there exists 
$\epsilon > 0$ such that $\zeta(s,\rho)$ has no
zero or  poles in 
$\mathrm{Re}(s) >  h(\rho) - \epsilon$  other than $h(\rho)$. 
Moreover, 
there exists $\alpha > 0$ such that 
we can bound 
$\log |\zeta(s,\rho)| = O(|\mathrm{Im}(s)|^\alpha)$ for 
$s$ satisfying $h(\rho) - \epsilon < \mathrm{Re}(s) < h(\rho)$ 
and $|\mathrm{Im}(s)| \geq 1$.
\end{thm}

\begin{proof}
That $\zeta(s,\rho)$ has a non-zero meromorphic extension to a half-plane
$\mathrm{Re}(s) > h(\rho) -\epsilon$, for some $\epsilon>0$, 
is a standard result (Corollary 10.6 of \cite{PP}).
The pole free strip is a consequence of the method of Dolgopyat \cite{dolgopyat} and the fact that,
as discussed above,
the limit set 
$\Lambda \subset \mathbb R P^{d-1}$ is a $C^1$ curve \cite{PSW}, \cite{ZZ}.
This method gives bounds on the iterates of the transfer operators $\mathcal L_s$ regarded as operators
 $\mathcal L_s: C^1\left(\coprod_{i=1}^k P_i, \mathbb C\right) 
 \to C^1\left(\coprod_{i=1}^k P_i, \mathbb C\right)$.  
 It is convenient to use the norm 
$$\|h\|_{1,t} = 
\begin{cases}
\max
 \left\{
\|h\|_\infty, 
\frac{\|h'\|_\infty}{|t|}
\right\} &\hbox{ if } |t| \geq 1\\
\max
 \left\{
\|h\|_\infty, 
\|h'\|_\infty
\right\} &\hbox{ if } |t| <  1
\end{cases}
$$ on $C^1\left(\coprod_{i=1}^k P_i, \mathbb C\right)$.
Following Proposition 7.4 of \cite{AGY}, 
the key ingredients required for the proof are that:
\begin{enumerate}
\item[(i)] $r: \coprod_{i=1}^k P_i \to \mathbb R$ is a $C^1$ function; and 
\item[(ii)] 
$r: \coprod_{i=1}^k P_i \to \mathbb R$ is not cohomologous to a constant, i.e. there is 
no $C^1$ function 
$u: \coprod_i P_i \to \mathbb R$
and constant $c \in \mathbb R$ such that $r = u \circ T - u + c$.
\end{enumerate}
We see that (i) holds by Lemma \ref{stuff} and
Lemma \ref{lem:limit_set_is_C1} and (ii) follows from Lemma \ref{key}
and Lemma \ref{mixingcondition}.

We can apply results on bounding iterates of transfer
operators: Corollary 2 and Corollary 3 of  \cite{dolgopyat}
and Proposition 4 of
\cite{PS-ajm},
to show that there exist constants $\sigma_0 < h$, $C > 0$ and $0 < \beta < 1$ such that whenever $s= \sigma + i t$ and 
$n = p [\log |t|] + l$, where $p \geq 0$ and 
$0 \leq l \leq [\log|t|]-1$, then 
$$
\|\mathcal L_{-sr}^n\|_{1,t} \leq 
C \beta^{p[\log |t|]} e^{l P(-\sigma r)}
$$
with respect to the norm $\|h\|_{1,t}
$.

By Lemma 2 of \cite{PS-ajm} (which, apart from the dependence on $|t|$ in the bound, appears in
\cite{ruelle-ihes}) we have the estimate
that for any 
$x_j \in P_j$, $j=1, \ldots, k$, and any $\beta_0$ satisfying
$\max\{\beta, \max_j\|D\psi_j\|_\infty\}< \beta_0 < 1$, there exists $C > 0$ such that 
$$
\left|
\sum_{|\underline i| =n} \det(D\psi_{\underline i}(z_{\underline i}))^s - \sum_{j=1}^k 
{\mathcal L}_s \chi_{P_j}(x_j)
\right| \leq C |t| n \beta_0^n, \quad \forall n \geq 1,
$$
where $\chi_{P_j}$ is the indicator function for $P_j$ and 
$\underline i$, $\psi_{\underline i}$ and $z_{\underline i}$ are as in Lemma \ref{traces}.
We can use this estimate to get a  bound on the logarithm of the absolute value of the zeta function of the form 
$
\log |\zeta(s, \rho)|= O\left( |t|^\alpha\right),
$
for some  $\alpha > 0$, as in (2.3) of \cite{PS-ajm}.
\end{proof}

In light of (1.4),  
and the observation that $Z(s+1, \rho)$ is uniformly bounded on the half-plane $\mathrm{Re}(s) > h(\rho) - \epsilon$, provided $\epsilon < 1$,
this gives  the following
generalization of 
Proposition
\ref{riemannhypothesisselberg}.

\begin{thm} 
Let $\rho : \Gamma \to \mathrm{SL}(d,\mathbb R)$ be a $(1,2)$-Anosov representation
which is $(1,1,2)$-hyperconvex.
Then there
exists 
$\epsilon > 0$ such that $Z(s,\rho)$ has no zero or poles in 
$\mathrm{Re}(s) >  h(\rho) - \epsilon$  other than the zero at $s= h(\rho)$. 
Moreover, there exists $\alpha > 0$ such that we can bound 
$\log |Z(s,\rho)| = O(|\mathrm{Im}(s)|^\alpha)$ for 
$s$ 
satisfying $h(\rho) - \epsilon < \mathrm{Re}(s) < h(\rho)$ 
and $|\mathrm{Im}(s)| \geq 1$.
\end{thm}

Although Theorem \ref{riemannhypothesis} is weaker than the full result known for the Selberg zeta function in the case of hyperbolic surfaces it nonetheless has analogous consequences for error terms in counting functions. To conclude the section, we address how the positions of the zeros influence asymptotic formulae.

\begin{definition}\label{counting}
For $T>0$, we define 
$$
\pi_\rho(T) = \#\{[g] \in \mathcal P  \hbox{ : } d_\rho([g]) \leq T \},
$$
the number of primitive conjugacy classes with weight at most $T$.
\end{definition}

It was shown in \cite{samCMH} (following the approach in
Chapter 6 of \cite{PP})
 that $\pi_\rho(T)$ satisfies the asymptotic formula
$$
\pi_\rho(T) \sim \frac{e^{h(\rho) T}}{h(\rho)T}, \quad \hbox{ as } T \to +\infty. 
$$
Using Theorem \ref{riemannhypothesis},
we can deduce a stronger result
with an error term provided we replace the leading term with a logarithmic integral.
As usual, we will write
\[
\mathrm{li}(x) := \int_2^x \frac{1}{\log u} du \sim \frac{x}{\log x}, \quad \mathrm{as}\  x \to +\infty.
\]
We then have the following result.

\begin{thm}\label{asymp}
There exists 
$\epsilon > 0$ such that
$$
\pi_\rho(T) = 
\mathrm{li}(e^{h(\rho)T})
\left( 1 + O(e^{-\epsilon T})\right).
$$
\end{thm}

\begin{proof}
Given Theorem \ref{riemannhypothesis}, we can apply the classical proof from number theory
(see Proposition 6 and Proposition 7 of \cite{PS-ajm}).
\end{proof}

\begin{remark}
It is worth noting that the trick used in Remark \ref{rem:adjoint_rep}
does not apply in the context of this section, since the adjoint representation of a $(1,1,2)$-hyperconvex
representation is not necessarily itself hyperconvex.
\end{remark}

\section{The $L$-function in higher rank}\label{sec:Lfunction}

We can generalize the definition of the zeta functions  to $L$-functions 
by incorporating additional  information associated to a unitary representation 
$U: \Gamma \to \mathrm{U}(N)$.  More precisely, we make the following definition.

\begin{definition}
We can associate to each projective Anosov representation
 $\rho: \Gamma \to  \mathrm{SL}(d, R)$ 
and to each unitary representation
 $R_\chi: \Gamma \to \mathrm{U}(N)$ (with character $\chi =\mathrm{Trace}(R_\chi)$),
 a Selberg type $L$-function formally defined by  
$$
Z(s,\rho,\chi)
= \prod_{n=0}^\infty \prod_{[g] \in \mathcal P} 
\det \left(
1- e^{-(s+n) d_{\rho}(g)} R_\chi([g])
\right), \quad s \in \mathbb C,
$$ 
and a Ruelle type $L$-function formally defined by
\[
\zeta(s,\rho,\chi)
= \prod_{[g] \in \mathcal P} 
\det \left(
1- e^{-s d_{\rho}(g)} R_\chi([g])
\right)^{-1}, \quad s \in \mathbb C,
\]
 where the products converge.
\end{definition}

If $R_{\chi_0}$ is the trivial representation then we see that these $L$-functions reduce to the 
corresponding zeta functions, i.e.
$Z(s, \rho, \chi_0) = Z(s,\rho)$ and $\zeta(s,\rho,\chi_0) = \zeta(s,\rho)$.
In the particular case of Fuchsian representations then these definitions reduce to the familiar definitions of 
$L$-functions for Fuchsian groups \cite{fried}.

It is easy to see that 
$Z(s,\rho,\chi)$ and $\zeta(s,\rho,\chi)$
converge for 
$\mathrm{Re}(s) > h({\rho})$.   We can also 
observe directly from the definitions that
$$
\zeta(s, \rho,\chi) = 
\frac{Z(s+1, \rho, \chi)}{Z(s, \rho, \chi)}. \eqno(5.1)
$$

The following result
generalizes Theorem   \ref{cor-main} and 
has a similar proof.

\begin{thm}\label{ruelleL}
Let  $\rho: \Gamma \to  \mathrm{SL}(d, \mathbb R)$ be a projective Anosov representation 
and let $R_\chi: \Gamma \to \mathrm{U}(N)$
be a unitary representation.
The Ruelle-type $L$-function 
$\zeta(s,\rho,\chi)$
 has the following properties.
 \begin{enumerate}
\item 
 $\zeta(s,\rho,\chi)$
 has a meromorphic extension to 
 the entire complex plane $\mathbb C$; and
  \item
  $\zeta(s,\rho,\chi)$
   has a simple  pole at $s=h(\rho)$ if and only if $R_\chi$ is trivial.
  \end{enumerate}
\end{thm}

The following result
generalizes Theorem   \ref{main} and 
has a similar proof.

\begin{thm}\label{mainLL}
Let  $\rho: \Gamma \to  \mathrm{SL}(d, \mathbb R)$ be a projective Anosov representation 
and let $R_\chi: \Gamma \to \mathrm{U}(N)$
be a unitary representation.
The generalized $L$-function $Z(s, \rho, \chi)$
 has the following properties.
 \begin{enumerate}
\item 
$Z(s,\rho,\chi)$
 has a meromorphic extension to 
 the entire complex plane $\mathbb C$; and
  \item
  $Z(s,\rho,\chi)$
   has a simple zero at $s=h(\rho)$ if and only if $R_\chi$ is trivial.
  \end{enumerate}
\end{thm}

 \end{document}